\documentclass[12pt,reqno,twoside]{amsart}

\usepackage[T1]{fontenc}
\usepackage{lmodern}
\usepackage{amsmath,amssymb,amsthm,mathtools}
\usepackage{microtype}
\usepackage{xcolor}
\usepackage[colorlinks,linkcolor=teal,citecolor=teal]{hyperref}

\newcommand{\C}{\mathbb C}
\newcommand{\R}{\mathbb R}
\newcommand{\cF}{\mathcal F}
\newcommand{\cV}{\mathcal V}
\newcommand{\cO}{\mathcal O}
\newcommand{\dbar}{\bar\partial}
\newcommand{\dd}{\mathrm{d}}
\newcommand{\st}{\mathrm{st}}
\newcommand{\loc}{\mathrm{loc}}
\newcommand{\poly}{\mathrm{poly}}
\newcommand{\Vol}{\operatorname{Vol}}

\newcommand{\fp}{\mathfrak{p}}

\newtheorem{theorem}{Theorem}
\newtheorem{lemma}{Lemma}
\newtheorem{proposition}{Proposition}
\newtheorem{problem}{Problem}
\theoremstyle{remark}

\numberwithin{equation}{section}

\title[Global nonsolvability of the leafwise $\dbar$-equation]{Global nonsolvability of the leafwise $\dbar$-equation on foliations with Stein leaves}
\author{Qingchun Ji}
\address{School of Mathematical Sciences, Fudan University, Shanghai 200433, China}
\email{qingchunji@fudan.edu.cn}
\author{Jun Yao}
\address{School of Mathematical Sciences, University of Electronic Science and Technology of China, Chengdu 611731, China}
\email{junyao@uestc.edu.cn}

\date{}

\begin{document}
	
\begin{abstract}
Motivated by our previous work on $q$-convexity and global solvability for Levi flat structures, we study the complex foliation associated with the non-tame family of Stein structures constructed by Forstneri\v c and Sigur\dh ard\'ottir. Although every leaf is closed and biholomorphic to $\C^2$ and the foliation admits a complete bundle-like Riemannian metric, there exists a smooth leafwise $\dbar$-closed $(0,1)$-form for which the leafwise $\dbar$-equation admits no global $L^r_{\loc}$ solution for any $0<r\le\infty$. This result also suggests a modification of Questions~2.10.4--2.10.6 in \emph{Foliations: Dynamics, Geometry and Topology}.
\end{abstract}
\maketitle

\section{Introduction}

Let $M$ be a smooth manifold and let $\cV\subseteq\C{\rm T}M$ be a formally integrable complex subbundle. Such structures and their associated differential complexes have been studied from the perspectives of solvability, cohomology, and regularity (see \cite{T92,CH01,BCH08,CJ21,ADd23a,AFJR24,JYY26,JY23',KN26} and the references therein). A formally integrable structure $\cV$ is called \emph{Levi flat} if $\cV+\overline{\cV}$ is a formally integrable complex subbundle. In our earlier work \cite{JY26}, we constructed a class of differential complexes associated with Levi flat structures, providing a global realization of the Treves complex and, in particular, of certain leafwise differential complexes \cite[Propositions~4.1 and~4.4]{JY26}. We also introduced a notion of $q$-convexity adapted to $\cV$ and established global $L^2_{\loc}$-solvability under the assumption that $M$ admits a smooth $q$-convex exhaustion function \cite[Theorem~1.1]{JY26}.

Although the convexity condition concerns derivatives tangent to the underlying foliation $\cV+\overline{\cV}$, the exhaustion hypothesis in \cite[Theorem~1.1]{JY26} is imposed on the whole manifold. This naturally raises the question of whether global solvability still holds if each leaf admits its own convex exhaustion function.

An important setting for this question is that of Levi flat CR structures. By Nirenberg's complex Frobenius theorem \cite{N57}, there are local adapted coordinates $(z,t)\in\C^n\times\R^k$ in which
\begin{equation*}
\cV|_U=\operatorname{span}_\C\left\{\frac{\partial}{\partial\bar z_j}:1\le j\le n\right\}.
\end{equation*}
The changes of adapted coordinates have the form
\[
(z,t)\longmapsto\bigl(\phi(z,t),\psi(t)\bigr),
\]
where $\phi(\cdot,t)$ is holomorphic. Such a structure is also called a \emph{complex foliation} $\cF$ in the sense of \cite[Definition~2.10.1]{EKA14}, with $\cV={\rm T}^{0,1}\cF$. For related integrability and regularity results, we 
refer the reader to Gong's parameter-dependent Frobenius--Nirenberg theorem \cite[Theorem~1.3]{G20}.

For $p,q\ge0$, let $\mathcal A^{p,q}_{\cF}(M)$ denote the space of smooth leafwise forms of type $(p, q)$, i.e.,
\[
\mathcal A^{p,q}_{\cF}(M)=C^\infty(M,\Lambda^{p,q}{\rm T}^*\cF),
\]
where ${\rm T}^*\cF=({\rm T}\cF)^*$. Let
\[
\dbar_{\cF}:\mathcal A^{p,q}_{\cF}(M)\longrightarrow\mathcal A^{p,q+1}_{\cF}(M)
\]
be the Cauchy--Riemann operator along the leaves defined by the $(p,q+1)$-component of the leafwise exterior derivative. The operator is therefore independent of the choice of adapted charts. In such coordinates $(z,t)$, for $u\in C^\infty(M)$,
\[
\dbar_{\cF} u=\sum_{j=1}^n\frac{\partial u}{\partial\bar z_j}\,\dd\bar z_j.
\]
For each fixed $p$, the leafwise Dolbeault complex is
\begin{equation}\label{eq:leafwise-complex}
0\longrightarrow\mathcal A^{p,0}_{\cF}(M)\xrightarrow{\dbar_{\cF}}
\mathcal A^{p,1}_{\cF}(M)\xrightarrow{\dbar_{\cF}}\cdots
\xrightarrow{\dbar_{\cF}}\mathcal A^{p,n}_{\cF}(M)\longrightarrow0.
\end{equation}
Since $\dbar_{\cF}^{2}=0$, it defines the leafwise (foliated) Dolbeault cohomology
\begin{equation*}
{\rm H}^{p,q}_{\cF}(M)=\frac{\ker\{\dbar_{\cF}:\mathcal A^{p,q}_{\cF}(M)\to\mathcal A^{p,q+1}_{\cF}(M)\}}{\dbar_{\cF}\mathcal A^{p,q-1}_{\cF}(M)},\quad q\ge1.
\end{equation*}

In particular, for $p=0$, the complex \eqref{eq:leafwise-complex} is naturally identified with a special case of the differential complex constructed in \cite[Proposition~4.4 and \S\,4.2.2]{JY26}. For this complex, a global convexity condition involving the leafwise complex Hessian ensures the vanishing of the corresponding $L^2_{\loc}$ leafwise cohomology in bidegrees $(0,q')$ with $q'\ge q$ \cite[(2.19), Definition~2.3, and Corollary~4.6]{JY26}.

Complex foliations with Stein leaves provide a natural framework for investigating global solvability of the leafwise $\dbar$-equation under leafwise convexity assumptions. Indeed, every Stein manifold admits a smooth strictly plurisubharmonic exhaustion function, and its Dolbeault cohomology vanishes in bidegrees $(p,q)$ with $q\ge1$ by the Dolbeault isomorphism and Cartan's Theorem B. However, solutions on the individual leaves need not combine into a solution on $M$ with the required transverse regularity or local integrability. El Kacimi Alaoui formulated this issue as Questions~2.10.4--2.10.6 of \cite[\S\,2.10]{EKA14}. These questions ask whether ${\rm H}^{0,q}_{\cF}(M)=0$ for $q\ge1$ in three different settings: foliations with closed Stein leaves; complete Riemannian foliations with closed Stein leaves; and differentiable products of a Stein manifold with an interval, with the leafwise complex structure allowed to vary. For explicit calculations of leafwise Dolbeault cohomology in concrete examples (see \cite{S08,EKS10,JJ23}).

In adapted coordinates, the leafwise $\dbar_{\cF}$-equation is a family of $\dbar$-equations depending on the transverse parameters. Diederich and Ohsawa studied the parameter dependence of canonical $L^2$ solutions \cite[Theorems~1 and~2]{DO91}. Gong and Kim established regularity results for parameter-dependent domains in $\C^n$ \cite[Theorem~1.1]{GK18}, while Forstneri\v c obtained parametric solvability for families of open Riemann surfaces \cite[Corollary~1.2]{F26}. For tame families of Stein structures, Forstneri\v c and Sigur\dh ard\'ottir proved global $\dbar$-solvability with continuous parameter dependence of all fibre derivatives \cite[Theorem~8.1]{FS26}.

We consider the family of Stein structures $J_t$ on $\C^2$ constructed by Forstneri\v c and Sigur\dh ard\'ottir \cite[Theorem~4.1]{FS26}, based on a non-Runge Fatou--Bieberbach domain of Wold \cite{W08}. Although every fibre is biholomorphic to $\C^2$, the holomorphic convex hulls of a fixed compact ball are not uniformly bounded as the parameter tends to zero. This phenomenon is used in \cite[Lemma~4.2]{FS26} to obtain restrictions on continuous families of global holomorphic functions. We exploit the same phenomenon to construct a smooth $\dbar_{\cF}$-closed form $f$ for which the equation $\dbar_{\cF}u=f$ admits no global solution in $L^r_{\loc}$ for any $0<r\le\infty$, with the equation understood in the sense of distributions on almost every leaf. More precisely, there exists a fixed point outside the closed unit ball that belongs to the $J_t$-holomorphic convex hull of that ball for all sufficiently small nonzero $t$. Combining this hull property with the submean inequality for holomorphic functions, we show that, for every $r>0$, any global $L^r_{\loc}$ solution $u$ would satisfy
\[
\|u(\cdot,t)\|_{L^r(B)}\geq c_r e^{1/t^2}\quad\text{for almost every sufficiently small }t>0,
\]
where $B\Subset\C^2$ is a fixed ball and $c_r>0$ is independent of $t$. This contradicts the local $L^r$-integrability of $u$ on the total space.

Our main global nonsolvability result for the leafwise $\bar\partial$-equation is stated below.

\begin{theorem}\label{thm:main}
	Let $\{J_t\}_{t\in(-1,1)}$ be the family of complex structures on $\C^2$ constructed by Forstneri\v c and Sigur\dh ard\'ottir \cite[Theorem~4.1]{FS26}. Consider the associated complex foliation $(M,\cF)$, where
	\[
	M=\C^2\times(-1,1),\qquad
	\cF=\{X_t=\C^2\times\{t\}\mid -1<t<1\},
	\]
	and each leaf $X_t$ is endowed with the complex structure $J_t$. Then the following hold:
	\begin{enumerate}
		\item[$(i)$] Every leaf $X_t$ is closed in $M$ and biholomorphic to $\C^2$, hence Stein.
		\item[$(ii)$] $(M,\cF)$ admits a complete bundle-like Riemannian metric.
		\item[$(iii)$] There exists a smooth leafwise $\dbar_{\cF}$-closed $(0,1)$-form $f$ such that the equation $\dbar_{\cF}u=f$ admits no global solution $u\in \bigcup\limits_{r>0}L^r_{\loc}(M)$.
	\end{enumerate}
\end{theorem}

In particular, ${\rm H}^{0,1}_{\cF}(M)\ne0$. Thus Theorem~\ref{thm:main} gives negative answers to Questions~2.10.4--2.10.6 of \cite{EKA14}. Moreover, the above complex foliation, viewed as a Levi flat CR structure, cannot admit a smooth $1$-convex exhaustion function in the sense of \cite[Definition~2.3]{JY26}. Otherwise, \cite[Corollary~4.6]{JY26} would provide an $L^2_{\loc}$ solution to the equation $\dbar_{\cF}u=f$.

In our previous work, a positivity assumption on a basic line bundle also yielded global solvability with finite-order Sobolev regularity for sufficiently high tensor powers on compact manifolds \cite[Theorem~1.2]{JY26}. By contrast, the global convexity result for Levi flat structures \cite[Theorem~1.1]{JY26} guarantees $L^2_{\loc}$-solvability, but does not address smooth global solvability. These results, together with Theorem~\ref{thm:main}, motivate the following question concerning global smooth solvability under the global $q$-convexity assumption in the sense of \cite[Definition~2.3]{JY26}.

\begin{problem}
	Let $(M,\cF)$ be a complex foliation of complex leaf dimension $n$, and let $1\le q\le n$. Suppose that $M$ admits a smooth exhaustion function $\varphi$ such that, for every leaf $F$ and every $x\in F$, the restriction of $\sqrt{-1}\,\partial\bar\partial(\varphi|_F)$ to $T_x^{0,1}F\cap\ker(\dd\varphi_x)$ has at least
	\[
	\dim_{\C}\bigl(T_x^{0,1}F\cap\ker(\dd\varphi_x)\bigr)-q+1
	\]
	positive eigenvalues. Does every $\dbar_{\cF}$-closed form $f\in\mathcal A^{0,q'}_{\cF}(M)$ with $q'\ge q$ admit a smooth global solution $u\in\mathcal A^{0,q'-1}_{\cF}(M)$ to $\dbar_{\cF}u=f$?
\end{problem}

For families of domains in $\C^n$ admitting a uniform plurisubharmonic exhaustion function, Gong and Kim \cite[Theorem~1.1(ii)]{GK18} established global $\dbar$-solvability with smooth dependence on the parameters. Their result gives an affirmative answer to the corresponding smooth-solvability question for such families of domains.

\section{Geometric preliminaries}

\subsection{A family of Stein structures}

We recall the specific construction in the proof of \cite[Theorem~4.1]{FS26}. Let $B_r$ denote the open Euclidean ball of radius $r>0$ centered at $0\in\C^2$, and set
\[
K=\overline{B_1}\subseteq\C^2.
\]
For a compact set $E\subseteq\C^2$, denote by $\widehat E_{\poly}$ its polynomial hull. The construction gives a holomorphic embedding
\[
\Phi:\C^2\hookrightarrow\C^2
\]
satisfying
\[
0\in\widehat{\Phi(K)}_{\poly}\setminus\Phi(\C^2),
\]
and a smooth family of diffeomorphisms
\[
\Psi_t:\C^2\to\Psi_t(\C^2)\subseteq\C^2,\quad t\in(-1,1),
\]
such that $\Psi_t(\C^2)=\C^2$ for every $t\ne0$ and
\begin{equation}\label{eq:maps}
	\Psi_0=\Phi,\quad
	\Psi_t=\Phi\ \text{on a neighbourhood of }
	\overline{B_{1/|t|}}\ \text{for}\ t\ne0.
\end{equation}
Let $J_{\st}$ denote the standard complex structure on $\C^2$, and set
\begin{equation}\label{eq:complex}
	J_t=(\dd\Psi_t)^{-1}\circ J_{\st}\circ\dd\Psi_t.
\end{equation}
The structures $J_t$ are integrable and depend smoothly on $t$. Moreover, $J_0=J_{\st}$, and each $(\C^2,J_t)$ is biholomorphic to $(\C^2,J_{\st})$. Since $\Phi$ is holomorphic for $J_{\st}$, the local identity in \eqref{eq:maps} yields
\begin{equation}\label{eq:standard}
 J_t=J_{\st}\ \text{on a neighbourhood of }\overline{B_{1/|t|}}\ \text{for}\ t\ne0.
\end{equation}
Moreover, the same construction gives the following holomorphic hull identity for $t\ne0$, recorded as formula (4.1) in \cite{FS26}:
\begin{equation}\label{eq:hull}
 \widehat K_{J_t}=\Psi_t^{-1}\bigl(\widehat{\Phi(K)}_{\poly}\bigr),
\end{equation}
where $\widehat K_{J_t}$ denotes the hull with respect to all functions in $\cO(\C^2,J_t)$. By \cite[Lemma~4.2]{FS26}, the set $\widehat{\Phi(K)}_{\poly}\cap\Phi(\C^2)$ is not relatively compact in $\Phi(\C^2)$. Since $\Phi$ is a biholomorphism onto its image, it follows that
\begin{equation}\label{eq:unbounded}
 \Phi^{-1}\bigl(\widehat{\Phi(K)}_{\poly}\bigr)\ \text{is unbounded in }\C^2.
\end{equation}

\subsection{The associated complex Riemannian foliation}

Recall that a Riemannian metric $g$ on a foliated manifold $(M,\cF)$ is called \emph{bundle-like} if its induced metric $g^\perp$ on the normal bundle ${\rm T}M/{\rm T}\cF$ satisfies
\[
\mathcal L_X g^\perp=0\ \text{for every }X\in\Gamma(M,{\rm T}\cF).
\]
A foliation admitting such a metric is called a \emph{Riemannian foliation} (see \cite[\S\,2.2.3]{EKA14} and \cite{Mo88}). We now verify the geometric conditions of Theorem~\ref{thm:main}.

\begin{proposition}\label{prop:geometry}
	Let $M=\C^2\times(-1,1)$, and let $\cF$ be the foliation by the fibres $X_t=\C^2\times\{t\}$, endowed with the complex structures $J_t$ defined in \eqref{eq:complex}. Then $(M,\cF)$ is a smooth complex foliation, every leaf is closed in $M$ and biholomorphic to $\C^2$, and $M$ admits a complete bundle-like Riemannian metric.
\end{proposition}

\begin{proof}
	Consider the smooth map
	\[
	\Theta:M\longrightarrow\C^2\times(-1,1),
	\quad \Theta(z,t)=(\Psi_t(z),t).
	\]
	Since
	\[
	\dd\Theta=
	\begin{pmatrix}
		\mathrm{d}_z\Psi_t & \ast \\
		0 & 1
	\end{pmatrix},
	\]
	and $\dd_z\Psi_t$ is invertible, the map $\Theta$ is a local diffeomorphism. By \eqref{eq:complex}, $w:=\Psi_t(z)$ is holomorphic along each leaf, so $(w,t)$ form smooth adapted coordinates. Thus $(M,\cF)$ is a smooth complex foliation in the sense of \cite[Definition~2.10.1]{EKA14}. Note that every leaf is a fibre of the projection $M\to(-1,1)$ and is therefore closed in $M$. Moreover, each leaf $(X_t,J_t)$ is biholomorphic to $(\C^2,J_{\st})$.

	Let $g_0$ denote the Euclidean metric on the underlying $\R^4$, and define a smooth family of $J_t$-Hermitian metrics by
	\[
	h_t(v,w)=g_0(v,w)+g_0(J_tv,J_tw),\quad v,w\in {\rm T}\C^2.
	\]
	We take the leafwise and transverse directions to be orthogonal, and set
	\begin{equation}\label{eq:metric}
		g=h_t+\frac{\dd t\otimes \dd t}{(1-t^2)^2}.
	\end{equation}
	The transverse metric induced by \eqref{eq:metric} depends only on $t$, so $g$ is bundle-like. Moreover, $h_t\ge g_0$, so
	\[
	g\ge g_{\mathrm{ref}}:=g_0+\frac{\dd t\otimes \dd t}{(1-t^2)^2}.
	\]
	Since $g\ge g_{\mathrm{ref}}$, the smooth exhaustion function
	\[
	\varphi(z,t)=\sqrt{1+|z|^2}-\log(1-t^2)
	\]
	satisfies $|\nabla^g\varphi|_g^2\le 5$. Thus $g$ is complete.
\end{proof}

\section{Proof of Theorem~\ref{thm:main}}

We begin with a consequence of \eqref{eq:hull} and \eqref{eq:unbounded} concerning a fixed point in the holomorphic hulls.

\begin{lemma}\label{lem:hull}
There exist ${\fp}\in\C^2$, $j\in\{1,2\}$, and $\delta_0>0$ such that $|z_j({\fp})|>1$ and, for every $0<|t|<\delta_0$ and every $h\in\cO(\C^2,J_t)$,
\begin{equation}\label{eq:hull-bound}
 |h({\fp})|\le\sup_{z\in K}|h(z)|.
\end{equation}
\end{lemma}

\begin{proof}
By \eqref{eq:unbounded}, we may choose ${\fp}\in\Phi^{-1}\bigl(\widehat{\Phi(K)}_{\poly}\bigr)$ and $j\in\{1,2\}$ with $|z_j({\fp})|>1$. For all sufficiently small $t\ne0$ we have ${\fp}\in B_{1/|t|}$, \eqref{eq:maps} then gives $\Psi_t({\fp})=\Phi({\fp})\in\widehat{\Phi(K)}_{\poly}$. By \eqref{eq:hull}, the point ${\fp}\in\widehat K_{J_t}$ for all these parameters $t$. If we take $\delta_0>0$ sufficiently small, then the defining inequality for the holomorphic hull yields the estimate \eqref{eq:hull-bound}.
\end{proof}

\begin{proof}[Proof of Theorem~\ref{thm:main}]
Let the point $\fp$, the index $j$ and the number $\delta_0>0$ be as in Lemma~\ref{lem:hull}. We introduce a leafwise $(0,1)$-form by
\begin{equation}\label{eq:rhs}
	f(z,t)=
	\begin{cases}
		e^{1/t^2}(\dbar_{J_t}z_j)(z), & t\ne0,\\
		0, & t=0.
	\end{cases}
\end{equation}
We now verify that $f$ is smooth on $M$. Fix $a\in\C^2$ and choose $R>|a|+1$. For $|t|<R^{-1}$ and $|z-a|<1$, we have $|z|<R<|t|^{-1}$ whenever $t\ne0$. It follows from \eqref{eq:standard} that
\[
\dbar_{J_t}z_j=\dbar z_j=0.
\]
Consequently, $f\equiv0$ on a neighbourhood of every point $(a,0)$. In particular, $f$ extends smoothly across $t=0$ as a leafwise form of type $(0,1)$. Write $f_t=f(\cdot,t)$. By the integrability of $J_t$, we obtain
\[
 \dbar_{J_t}f_t=e^{1/t^2}\dbar_{J_t}^{\,2}z_j=0\ \text{for}\ t\ne0,
\]
and the same holds at $t=0$ by continuity. Hence $\dbar_{\cF}f=0$.

Let $0<r<\infty$. Assume that $u\in L^r_{\loc}(M)$ satisfies $\dbar_{\cF}u=f$ in the sense of distributions on almost every leaf. By Tonelli's theorem, $u_t:=u(\cdot,t)\in L^r_{\loc}(\C^2)$ for almost every $t\in(-1,1)$. Moreover,
\begin{equation}\label{eq:slice}
	\dbar_{J_t}u_t=f_t\quad\text{for almost every }t\in(-1,1).
\end{equation}
When $r\ge1$, since $u\in L^1_{\loc}(M)$ and $\dbar_{\cF}$ involves no transverse derivatives, the leafwise equations in \eqref{eq:slice} are equivalent to $\dbar_{\cF}u=f$ in the sense of distributions on $M$ by Fubini's theorem.

For such $t\ne0$, set
\[
F_t=z_j-e^{-1/t^2}u_t.
\]
By \eqref{eq:slice} and \eqref{eq:rhs}, we have
\[
\dbar_{J_t}F_t=0.
\]
By the $\dbar$-Weyl lemma in local $J_t$-holomorphic coordinates, $F_t$ has a $J_t$-holomorphic representative, which we still denote by $F_t$.

Choose $R>\max\{1,|{\fp}|\}+2$ and fix $0<\delta<\min\{1,\delta_0,(R+1)^{-1}\}$. For $0<t<\delta$, \eqref{eq:standard} and \eqref{eq:rhs} give
\[
J_t=J_{\st},\quad f_t=0\quad\text{on }B_{R+1}.
\]
Hence, for almost every $t\in(0,\delta)$, \eqref{eq:slice} yields
\[
\dbar u_t=0\ \text{on }B_{R+1}
\]
in the sense of distributions. By the $\dbar$-Weyl lemma, $u_t$ has a holomorphic representative on this ball. Since $F_t$ is also holomorphic there, the identity $F_t=z_j-e^{-1/t^2}u_t$ holds pointwise for their holomorphic representatives.

Applying the hull inequality \eqref{eq:hull-bound}, we obtain, for almost every $t\in(0,\delta)$,
\[
  |z_j({\fp})-e^{-1/t^2}u_t({\fp})|=|F_t({\fp})|\le\sup_K|F_t|
  \le 1+e^{-1/t^2}\sup_K|u_t|.
\]
The reverse triangle inequality therefore implies
\begin{equation}\label{eq:gap}
	|z_j({\fp})|-1\le e^{-1/t^2}\bigl(|u_t({\fp})|+\sup_K|u_t|\bigr)\quad \text{a.e. in}\ (0,\delta).
\end{equation}
Since $u_t$ is holomorphic on $B_{R+1}$ and
$K\cup\{\fp\}\subseteq B_R$, the submean inequality for $|u_t|^r$ gives
\[
|u_t(x)|\le
\Vol(B_{1/2})^{-1/r}
\|u_t\|_{L^r(B_{R+1})}
\quad\text{for every }x\in K\cup\{\fp\},
\]
where $\Vol$ denotes Euclidean volume. Combining this estimate with \eqref{eq:gap}, we obtain
\begin{equation}\label{eq:blowup}
	\|u_t\|_{L^r(B_{R+1})}\ge c_r e^{1/t^2}
	\quad\text{for almost every }t\in(0,\delta),
\end{equation}
where
\[
c_r=\frac{|z_j(\fp)|-1}{2}
\Vol(B_{1/2})^{1/r}>0
\]
is independent of $t$.

Since $\overline{B_{R+1}}\times[0,\delta]$ is compact in $M$, the local $L^r$-integrability of $u$ and Tonelli's theorem imply
\[
\int_0^\delta\|u_t\|_{L^r(B_{R+1})}^r\,\dd t<\infty.
\]
However, \eqref{eq:blowup} forces
\[
\int_0^\delta\|u_t\|_{L^r(B_{R+1})}^r\,\dd t
\ge c_r^r\int_0^\delta e^{r/t^2}\,\dd t=+\infty,
\]
which is a contradiction. Thus there is no such solution for any $0<r<\infty$. The case $r=\infty$ follows from $L^\infty_{\loc}(M)\subseteq L^1_{\loc}(M)$. Together with Proposition~\ref{prop:geometry}, this completes the proof of Theorem~\ref{thm:main}.
\end{proof}

\section*{Statements and Declarations}

\noindent\textbf{Funding.} The first author is supported by the National Natural Science Foundation of China under Grant Nos. 12431005 and 12571086; the second author is supported by the National Natural Science Foundation of China under Grant No. 12601141.

\medskip
\noindent\textbf{Competing Interests.} The authors declare that they have no competing interests.

%\medskip
%\noindent\textbf{Data Availability.} Data sharing is not applicable to this article as no datasets were generated or analyzed during the current study.

\end{document}